\documentclass[11pt, a4paper, reqno]{amsart}
\usepackage{amsthm,amssymb,amsmath}
\usepackage{graphicx}
\usepackage{enumerate}
\usepackage{hyperref}

\theoremstyle{definition}
\newtheorem{definition}{Definition}
\newtheorem{example}[definition]{Example}
\theoremstyle{remark}

\theoremstyle{plain}
\newtheorem{lemma}[definition]{Lemma}

\newtheorem{proposition}[definition]{Proposition}
\newtheorem{theorem}[definition]{Theorem}
\newtheorem{cor}[definition]{Corollary}
\newtheorem{conjecture}[definition]{Conjecture}

\newcommand{\set}[1]{\left\{{#1}\right\}}

\newcommand{\vek}[1]{\mathbf{#1}}

\newcommand\setsuchas[2]{\left\{\,{#1}\,\vrule\,{#2}\,\right\}}

\newcommand{\Nat}{{\mathbb{N}}}

\newcommand{\C}{{\mathbb{C}}}

\newcommand{\A}{\mathcal{A}}
\newcommand{\PB}{\mathcal{P}}
\newcommand{\HI}{\mathcal{H}}
\newcommand{\DO}{\mathcal{D}}
\newcommand{\tdeg}[1]{\lvert {#1} \rvert}
\newcommand{\del}{\partial}
\newcommand{\e}{\underline{\vek{e}}}
\newcommand{\sbinom}[2]{\left( \! \left( {\genfrac{}{}{0pt}{}{#1}{#2}} \right) \! \right)}
\newcommand{\hadamard}{\odot}
\newcommand{\tensor}{\otimes_S}
\newcommand{\Tor}{\mathrm{Tor}_S}
\begin{document}

\title[Differential operators on a generic hyperplane arrangement]{A
  conjecture on  Poincar{\'e}-Betti series of  modules of 
  differential operators on a generic hyperplane arrangement} 
\author{Jan Snellman}
\address{Department of Mathematics, Stockholm University\\
SE-10691 Stockholm, Sweden}
\email{Jan.Snellman@math.su.se}
\date{\today}
\begin{abstract}
  Holm \cite{Holm:Paper1, Holm:Diff} studied  modules of higher
  order differential operators (generalizing derivations) on generic
  (central) hyperplane arrangements. We use his results to determine
  the Hilbert series of these modules. We also give a conjecture about
  the Poincar{\'e}-Betti series; these are known for the module of
  derivations through the work of Yuzvinsky \cite{Yuz:FreeRes} and
  Rose and Terao \cite{RoseTerao}.
\end{abstract}
\maketitle

\sloppy

\begin{section}{Introduction}
  The module of derivations \(\DO^{(1)}(\A)\) of a hyperplane
  arrangement \(\A \in \C^n\) (henceforth called an  \(n\)-arrangement) is an 
  interesting and much studied object \cite{OrlikTerao:ArrHyp,
    GeomFreeArr, BasicDer}. In particular, the question whether this
  module is free, for various classes of arrangements, has received
  great attention.

  On the other hand, the module of higher differential operators
  \(\DO^{(m)}(\A)\) received their first incisive treatment in the PhD
  thesis of P{\"a}r Holm \cite{Holm:Diff}. The deepest result in that
  work concerns so-called \emph{generic arrangements}, which are
  arrangements where every intersection of \(s \le n\) hyperplanes in \(\A\)
  have the expected codimension \(s\). Holm gave a concrete generating
  set for \(\DO^{(m)}(\A)\), proved an extension of Saito's
  determinental criteria for freeness of derivations,  and used these
  results to tackle the question of 
  \emph{higher order freeness} for generic arrangements, i.e. the
  question when \(\DO^{(m)}(\A)\) is a free module (in which case we
  say that \(\A\) is \emph{\(m\)-free}). In brief, he
  showed that  
  \begin{enumerate}[(i)]
  \item   all 2-arrangements are \(m\)-free for all \(m\), 
  \item  all \(n\)-arrangements with \(\tdeg{A} \le n\) are \(m\)-free
    for  all \(m\),
  \item \label{case:interesting} 
    if \(n \ge 3\) and \(r > n\) and \(m < r-n+1\), then 
    \(\A\) is not \(m\)-free,
  \item if \(n \ge 3\) and \(r > n\) and \(m = r-n+1\), then 
    \(\A\) is  \(m\)-free.
  \end{enumerate}
  He conjectured that if \(n \ge 3\) and \(r > n\) and \(m > r-n+1\) then 
  \(\A\) is  \(m\)-free.
  
  For \(m=1\), \eqref{case:interesting} becomes \(r > n \ge 3\).
  Yuzvinsky \cite{Yuz:FreeRes} and independently Terao and Rose
  \cite{RoseTerao} have showed that the modules of derivations of
  generic arrangements are non-free, with a minimal free resolution of
  length \(n-2\). More precisely, they have showed that the graded
  Betti numbers are given by 
  \begin{displaymath}
    \beta_{k,u} = 
    \begin{cases}
      1 & \text{ if } k=0 \text{ and } u=1 \\
      \binom{r}{n-k} \binom{r-n+k-2}{k-1} & \text{ if } u+n-r-1=0
      \\
      0 & \text{ otherwise }
    \end{cases}
  \end{displaymath}
  so that the Poincar{\'e}-Betti series can be expressed as 
  \begin{displaymath}
b + \left[ t^{n-r-1}(1+bt)^{r-1}
          \right] \hadamard (1-t)^{n-r},
  \end{displaymath}
  where \(t\) enumerates homological degree, \(b\) ring degree, and
  \(\hadamard\) denotes Hadamard product of power series. 
  
  This articles contains a short expos{\'e} of various ways of
  calculating modules of differential operators (on a computer),  a
  brief review of the work of Holm,  and finally some conjectures
  supported by extensive computer experiments. The most important one
  is the conjectured formula for the Poincar{\'e}-Betti series of
  \(\DO^{(m)}(\A)\) when \(\A\) is a generic \(n\)-arrangement with
  \(\tdeg{\A}=r\), \(3 \le n\), \(r \ge m+n\):
  \begin{multline*}
        \PB(\DO(\A_{n,r})^{(m)}) = b^m + \\
        + \left\{ t^{-r +n - 1} \bigl((1+bt)^m -(bt)^m \bigr) 
        \bigr(1+bt\bigr)^{r-m} \right\} \hadamard (1-t)^{m - r+n - 1}.
  \end{multline*}

\end{section}

\begin{section}{Notation}
  For basic terminology regarding hyperplane arrangements, we refer to
  Orlik and Terao's treatise \cite{OrlikTerao:ArrHyp}. For more
  details on the Grothendieck ring of differential operators, see for
  instance the PhD thesis by Holm \cite{Holm:Diff}, or
  \cite{Bjork:DiffOp,Coutinho:Primer}.

  Let \(\A\) be an affine central hyperplane arrangement in \(\C^n\), with
  defining polynomial \(p=\prod_{i=1}^r p_i \in  S=\C[x_1,\dots,x_n]\). 
  We let \(\DO(S)\) denote the Weyl algebra of differential operators on
  \(S\). This is the set of all finite \(S\)-linear combinations
  \begin{equation}
    \label{eq:rm}
    \delta = \sum_{\alpha \in \Nat^n} c_\alpha \del^\alpha
  \end{equation}
  An element of \(\DO(S)\) can be regarded as a partial differential
  operator with polynomial coefficients, and thus it induces  an
  \(S\)-algebra endomorphism. We use 
  the notation \(Q * v\) to denote the action of \(Q \in \DO(S)\) on \(v
  \in S\). 

  \(\DO(S)\) is an \(S\)-module in a natural way: the action is given by 
  \begin{displaymath}
    q  \sum_{\alpha \in \Nat^n} c_\alpha \del^\alpha =
    \sum_{\alpha \in \Nat^n} q c_\alpha \del^\alpha.
  \end{displaymath}
  If  in \eqref{eq:rm} all \(\alpha\) with \(c_\alpha \neq 0\) have
  total degree \(m\)
  we say that \(\delta\) is a homogeneous \(m\)'th order operator, and
  write \(\tdeg{\delta}=m\), or \(\delta \in \DO^{(m)}(S)\).
  If, in addition, all \(c_\alpha\) occurring in \eqref{eq:rm} are homogeneous
  polynomials of total degree \(v\), we say that \(\delta\) is
  homogeneous of polynomial degree \(v\). Thus \(\DO(S)\) is a bi-graded
  \(S\) module, where we use the convention that \(x^\alpha
  \del^\beta\) has bi-grade  
  \((k,m)=(\tdeg{\alpha},\tdeg{\beta})\).

  Let \(m\) be a positive integer. We put 
  \begin{equation}
    \label{eq:drm}
    \begin{split}
    \DO(\A) &= \setsuchas{\delta = \sum_{\alpha \in \Nat^n} c_\alpha
      \del^\alpha}{c_\alpha \in S,\,\,
     \delta * \langle p \rangle \subseteq \langle p   \rangle 
    }  \\
    \DO^{(m)}(\A) &= \setsuchas{\delta = \sum_{\tdeg{\alpha}=m} c_\alpha
      \del^\alpha}{c_\alpha \in S,\,\,
     \delta * \langle p \rangle \subseteq \langle p   \rangle 
    } 
    \end{split}
  \end{equation}

  In particular, \(\DO^{(1)}(\A)\) is the much-studied module of
  derivations of \(\A\).

  It is a fact that \(\DO^{(m)}(\A)\) is a graded
  \(S\)-module, where the \(\Nat\)-grading is given by polynomial
  degree. Holm \cite{Holm:Paper1} 
  showed that \(\DO(\A) = \bigoplus_{m \ge 0} \DO^{(m)}(\A)\).
  Hence, the \(S\)-module
  \(\DO(\A)\) is bi-graded, and hence so is all Tor modules.
  The \(S\)-module \(\DO(\A)\) is not necessarily finitely generated, but
  every \(\DO^{(m)}(\A)\) is. Consequently we can calculate the graded
  minimal free resolution 
  \begin{equation}
    \label{eq:freeres}
    0 \leftarrow \DO^{(m)}(\A) \leftarrow 
    \bigoplus_{i} \beta_{1,i}^{(m)} S(-i) 
    \leftarrow  \cdots \leftarrow
    \bigoplus_{i} \beta_{\ell,i}^{(m)} S(-i) \leftarrow 0
  \end{equation}
  where, by the Hilbert syzygy theorem, \(\ell \le n\).
  We define the Poincar{\'e}-Betti series and Hilbert series of
  \(\DO^{(m)}(\A)\) and of \(\DO(\A)\) by 
  \begin{equation}
    \label{eq:pb}
    \begin{split}
      \PB(\DO^{(m)}(\A))(b,t) &= \sum_{j,i}  \beta_{j,i}^{(m)} b^it^j \\
      \PB(\DO(\A))(a,b,t) &= \sum_{m} a^m\PB(\DO^{(m)}(\A))(b,t) \\
      \HI(\DO(\A))(a,b) &= (1-b)^{-n}\PB(\DO(\A))(a,b,-1) \\
      \HI(\DO^{(m)}(\A))(b) &= (1-b)^{-n} \PB(\DO^{(m)}(\A))(b,-1)
    \end{split}
  \end{equation}


  \begin{subsection}{Additional notation}
    We define \(\sbinom{a}{b} = \binom{a+b-1}{b}\),
    i.e. \(\sbinom{a}{b}\) is the number of multisets of weight \(b\) on
    an \(a\)-set.    

    If \(f(t) = \sum_{i=0}^\infty c_i t^i\) is a formal power series
    in \(t\), with coefficients in some commutative ring, we put
    \([t^\ell] f(t) = c_\ell\).

    \begin{definition}
      Let \(R\) be a commutative ring and let \(f, g \in
      R[[t^{-1},t]]\) be two formal Laurent series, i.e.
      \begin{displaymath}
        f = \left( \sum_{i=-\infty}^\infty a_it^i \right),  \qquad
        g = \left( \sum_{i=-\infty}^\infty b_it^i \right).
      \end{displaymath}
      We define the \emph{Hadamard product} of \(f\) and \(g\) by 
      \begin{displaymath}
        f \hadamard g =  \sum_{i=-\infty}^\infty a_ib_it^i.  
      \end{displaymath}
    \end{definition}

  \end{subsection}

\end{section}

\begin{section}{Calculating modules of differential operators on
    hyperplane arrangements}
  We shall review some methods of calculating generators of the
  \(S\)-module \(\DO^{(m)}(\A)\).

  \begin{subsection}{The ``Jacobian'' method}
  First, we describe the most straight-forward method,
  a slight variant of which is used in the Macaulay 2
  package ``D-modules.m2'' \cite{Dmod} by Harry Tsai and Anton Leykin.

  \begin{lemma}
    Let \(\DO^{(m)}(S) \ni \delta =  \sum_{\tdeg{\alpha}=m} c_\alpha
    \del^\alpha\) be homogeneous of polynomial degree \(v\). Then 
    \(\delta \in \DO^{(m)}(\A)\) iff \(\delta * (x^\beta p) \in \langle p
    \rangle\) for all \(\tdeg{\beta} < v\).  
  \end{lemma}

  \begin{cor}
    Let \(m\) be a positive integer. Let \(G\) be a row matrix 
    whose entries are  \(\setsuchas{\del^{\alpha}}{\tdeg{\alpha}=m}\),
    let \(H\) be a column matrix whose  entries are
    \(\setsuchas{x^\beta}{\tdeg{\beta}  < m}\),
    and let \(A\) be the matrix indexed by \(G\) and \(H\) where the
    \((\del^\alpha,x^\beta)\) entry is \(\del^\alpha * (x^\beta p)\).
    Let \(B = [A | pI]\), where \(I\) is the identity matrix of
    appropriate dimension. Then the syzygy module of the columns of \(B\)
    correspond to \(\DO^{(m)}(\A)\). More
    precisely, if 
    \begin{displaymath}
      [A | pI]  \begin{bmatrix}\vek{u} \\ \vek{w} \end{bmatrix} = \vek{0}
    \end{displaymath}
    then 
    \begin{displaymath}
      \sum_{\alpha} u_\alpha \del_\alpha \in \DO^{(m)}(\A),
    \end{displaymath}
    and this is an isomorphism of \(S\)-modules.
  \end{cor}

  \begin{example}\label{ex:d2syz}
    Suppose that \(S=\C[x_1,x_2]\) and that \(p=x_1\). For \(m\)=2
    the matrix \(A\) is 
    \begin{center}
    \begin{tabular}{c|cccc}
      & \(\del_1^2\) & \(\del_1\del_2\) & \(\del_2^2\) \\
      \hline
      1 &  0 &0& 0 \\
      \(x_1\)&  2 &0 &0\\
      \(x_2\)& 0 & 1 & 0
    \end{tabular}
    \end{center}
    so in order to calculate \(\DO^{(2)}(\A)\) we should calculate the
    syzygies of the columns of 
    \begin{displaymath}
      B = 
      \begin{pmatrix}
       0 &0& 0 & x_1 & 0 & 0\\
       2 &0 &0 & 0 & x_1 & 0\\
       0 & 1 & 0 & 0 & 0 & x_1
      \end{pmatrix}
    \end{displaymath}
    A generating set is 
    \begin{displaymath}
      \begin{bmatrix}
        0 \\       0 \\   1\\ 0\\ 0\\ 0 
      \end{bmatrix},
      \begin{bmatrix}
       -1/2x_1 \\0 \\   0\\ 0\\ 1\\ 0
      \end{bmatrix},
      \begin{bmatrix}
       0 \\ -x_1\\ 0\\ 0\\ 0\\ 1
      \end{bmatrix},
    \end{displaymath}
    so \(\DO^{(2)}(\A)\) is generated by \(\del_2^2\), \(-1/2 x_1\del_1^2\),
    and \(-x_1 \del_1\del_2\).
  \end{example}
    
  \end{subsection}

  \begin{subsection}{The method of intersecting modules}
    The following result is proved in  \cite{Holm:Paper1}.
    \begin{theorem}
      \(\DO^{(m)}(p) = \cap_{i=1}^r \DO^{(m)}(p_i)\).
    \end{theorem}

    Holm also showed how to calculate \(\DO^{m}(p_i)\). 
    Let \(q \in S_1\) be a linear form, and let \(H \subset \C^n\)
    be its associated linear variety (a hyperplane through the
    origin). 
    \begin{definition}\label{def:V}
    Let \(V\) be the \(\C\)-vector space \(V= \sum_{i=1}^n \C
    \del^i\) and define 
    \begin{displaymath}
      V_H = \setsuchas{\sum_{i=1}b_i \del_i}{(b_1,\dots,b_n) \in H}.
    \end{displaymath}
    Then \(V_H\) is a codimension 1 subspace of \(V\).
    \end{definition}
    
    \begin{lemma}
      Let \(\mathcal{N}\) be the module of derivations annihilating
      \(q\). Then \(\mathcal{N} = S V_H\), and if \(\delta\) is any
      derivation such that \(\delta * q = a q\), \(a \in \C^*\), then
      \(\DO^{(1)}(q) = \mathcal{N} + S\delta\).
    \end{lemma}

    \begin{proposition}\label{prop:singlehyp}
      Let \(M=\set{\delta_1,\dots, \delta_n}\) be a basis for \(V\) such that
      \(\set{\delta_1,\dots,\delta_{n-1}}\) is  a basis for
      \(V_h\). Then \(N=\set{\delta_1,\dots,\delta_{n-1},q\delta_n}\)
      generates \(\DO^{(1)}(q)\), and
      \begin{equation}
        \label{eq:dmgen}
        \DO^{(m)}(q) = 
        \sum_{\substack{\tdeg{\alpha}=m \\ \alpha_n=0}} S
        \delta^\alpha + 
        \sum_{\substack{\tdeg{\alpha}=m \\ \alpha_n>0}} Sq  \delta^\alpha 
      \end{equation}
    \end{proposition}

    The above results can be succinctly summarized as follows: 
    let \(A(p)\) be the \(n \times n\) coordinate matrix of
    \(N\), i.e. the matrix formed by the
    coordinate vectors of elements of \(N\), and let \(S^m A(p)\)
    denote the \(m\)'th symmetric power. Let \(B^m(p)\) be the result
    of replacing any     occurrence of \(q^i\), with \(i > 1\), by
    \(q\). In other words, if \(A(p)\) is regarded as the matrix of
    an endomorphism \(\phi: S^n \to S^n\), then \(S^m A(p)\) is the
    matrix of the endomorphism \(S^m \phi: S^m S^n \to S^m S^n\), and 
    \(B^m(p)\) is the matrix of the associated endomorphism on \(S^m
    T\), 
where \(T=S/(q-q^2)\).
    
    \begin{example}[Example \ref{ex:d2syz} cont.]
      If \(S=\C[x_1,x_2]\) and  \(p=x_1\), then \(V_H\) is spanned by
      \(\del_2\).
      Hence, we can take \(M=\set{\del_2,\del_1}\) and
      \(N=\set{\del_2,x_1\del_1}\), so if we order the monomials of
      degree two as \(x_1^2,x_1x_2,x_2\) then
      \begin{displaymath}
        A = 
        \begin{bmatrix}
          0 & x_1 \\ 1 & 0
        \end{bmatrix},
        \quad
        S^2 A = 
        \begin{bmatrix}
          0 & 0   & x_1^2  \\
          0 & x_1 & 0 \\
          1 & 0  & 0
        \end{bmatrix},
        \quad
        B = 
        \begin{bmatrix}
          0 & 0   & x_1  \\
          0 & x_1 & 0 \\
          1 & 0  & 0
        \end{bmatrix}.
      \end{displaymath}
      Thus we recover the result that \(\DO^{(2)}(x_1)\) is generated by
      \begin{displaymath}
        [\del^2, \del_1\del_2, \del_2^2] B = [\del_2^2, x_1
        \del_1\del_2, x_1 \del_1^2].
      \end{displaymath}
    \end{example}

    Now recall (see for instance \cite[Theorem 3.8.3]{LouAD:GB})
    the following method of computing
    intersection of submodules of free modules. Suppose that
    \(M_1,\dots,M_\ell\) are 
    submodules of the free module \(S^s\), that
    \(\e=[\vek{e}_1,\dots,\vek{e}_s]\) is a basis of \(S^s\), and that
    the matrix \(A_i\) consists of the coordinate vectors (as column
    vectors) for a generating set of \(M_i\). Then the truncations of the
    syzygies of the matrix 
    \begin{displaymath}
      \begin{bmatrix}
        I_s & A_1 & 0 & 0 & \cdots & 0 \\
        I_s & 0 & A_2 & 0& \cdots & 0 \\
        \vdots & \vdots & \vdots &  \ddots & \cdots &\vdots \\
        I_s & 0 & 0 & 0& \cdots & A_s \\
      \end{bmatrix}
    \end{displaymath}
    correspond to elements in \(\cap_{i=1}^s M_i\).

    \begin{example}
      We have that \(\DO^{(2)}(x_1x_2) = \DO^{(2)}(x_1) \cap
      \DO^{(2)}(x_2)\).
      With respect to the basis \([\del_1, \del_1\del_2, \del_2^2]\)
      for \(\DO^{(2)}(S)\) the matrices of \(\DO^{(2)}(x_1)\) and
      \(\DO^{(2)}(x_2)\) can be taken to be
      \begin{displaymath}
        \begin{bmatrix}
          x_1 & 0 & 0 \\
          0 & x_1 & 0 \\
          0 & 0 & 1
        \end{bmatrix}
        \text{ and }
        \begin{bmatrix}
          1 & 0 & 0 \\
          0 & x_2 & 0 \\
          0 & 0 & x_2
        \end{bmatrix}.
      \end{displaymath}
      The syzygies of the matrix 
      \begin{displaymath}
        \begin{bmatrix}
1 & 0 & 0 &          x_1 & 0   & 0 & 0 & 0 & 0\\
0 & 1 & 0 &         0   & x_1 & 0 & 0 & 0 & 0\\
0 & 0 & 1 &          0   & 0   & 1 & 0 & 0 & 0\\
1 & 0 & 0 &          0   & 0   & 0 & 1 & 0 & 0 \\
0 & 1 & 0 &          0   & 0   & 0 & 0 & x_2 & 0 \\
0 & 0 & 1 &          0   & 0   & 0 & 0 & 0 & x_2
        \end{bmatrix}
      \end{displaymath}
      are generated by 
      \begin{displaymath}
        \begin{bmatrix}
         0 \\  0 \\ -x_2 \\ 0\\ 0\\  x_2\\ 0\\ 0\\  1 
        \end{bmatrix}
        , \,
        \begin{bmatrix}
           -x_1 \\ 0\\  0\\  1\\ 0\\  0\\ x_1\\ 0\\  0 
        \end{bmatrix}
        , \,
        \begin{bmatrix}
            0 \\  x_1x_2\\   0\\  0\\ -x_2\\ 0\\ 0\\ -x_1\\ 0
        \end{bmatrix},
      \end{displaymath}
      so the intersection of the two modules is generated by 
      \begin{displaymath}
        \begin{bmatrix}
         0 \\  0 \\ -x_2 
        \end{bmatrix}
        , \,
        \begin{bmatrix}
           -x_1 \\ 0\\  0
        \end{bmatrix}
        , \,
        \begin{bmatrix}
            0 \\  x_1x_2\\   0
        \end{bmatrix},
      \end{displaymath}
      and \(\DO^{(2)}(x_1x_2)\) is generated by 
      \begin{math}
-x_2 \del_2,        -x_1\del_1,  x_1x_2 \del_1\del_2.
      \end{math}
    \end{example}
  \end{subsection}

\begin{subsection}{Calculating modules of differential operators on
    generic arrangements} 
    Holm \cite{Holm:Paper1} gives the following method for
    constructing a (not necessarily minimal) generating set of
    \(\DO^{(m)}(\A)\), when \(\A\) is a generic \(n\)-arrangement.

    \begin{subsubsection}{The case \(r > n\)}
    Holm  showed that for any positive
    \(m\), and for any central arrangement \(\A\), the modified Euler
    derivation  
    \begin{displaymath}
      \varepsilon_m = \sum_{\tdeg{\alpha}=m} \frac{m!}{\alpha!}
      x^\alpha \del^\alpha 
    \end{displaymath}
    belongs to \(\DO^{(m)}(\A)\). He then proceeded to  find other
    generators as follows.

    Recall the
    definition of \(V\) and of \(V_H\) from definition~\ref{def:V}. We
    let \(H_i\) be the hyperplane associated with the linear form
    \(p_i\), and put \(V_i = V_{H_i}\). Choose a basis element for
    each intersection of \(n-1\) of the \(V_i\)'s (by genericity, this
    intersection is 1-dimensional) and let
    \(M=\set{\delta_1,\dots,\delta_t}\), with \(t=\binom{r}{n-1}\), be
    the set of all of these. We define a subset \(D\) of the
    derivations on \(S\) by 
    \begin{displaymath}
      D=\set{P_1\delta_1,\dots,P_t\delta_t},
    \end{displaymath}
    where each \(P_i\) is the product of those \(p_j\):s which are not
    annihilated by \(\delta_i\). If all \(p_j\):s  are 
    annihilated by \(\delta_i\), we put \(P_i=1\).

    Holm \cite{Holm:Paper1} showed that \(D \subset \DO^{(1)}(\A)\),
    and furthermore that:
    \begin{theorem}\label{thm:gengen}
      Suppose \(P= p_1\cdots p_r \) is the defining polynomial of a generic
      arrangement \(\A\). Let \(M\) and \(D\) be as above, and let
      \(I\) be the principal ideal on \(P\). Then 
      \begin{displaymath}
        \DO^{(m)}(\A) = \bigoplus_{m \ge 0} \left( \sum_{\tdeg{\alpha}=m}
          (I:(I:P^\alpha))\delta^\alpha + S\varepsilon_m \right)
      \end{displaymath}
      as a bi-graded \(S\)-module.
    \end{theorem}
    Here, the double colon ideal \((I:(I:P^\alpha))\) is a principal
    ideal in \(S\), with the generator given by \(p_{i_1} \cdots
    p_{i_\ell}\), the product of those \(p_j\) such that some
    \(\delta_i\) with \(\alpha_i \neq 0\) does not annihilate \(p_j\).

    \begin{example}
      Let \(p=xy(y-x) \in \C[x,y]=S\). Then the associated arrangement
      is generic. We get that 
      \begin{displaymath}
        P_1= y(y-x), \quad P_2=x(y-x), \quad P_3=xy, 
      \end{displaymath}
      and 
      \begin{displaymath}
        \delta_1=\del_y, \quad \delta_2=\del_x, \quad \delta_3=\del_x+\del_y.
      \end{displaymath}
      Hence \(\DO^{(0)}(I)=S\), and \(\DO^{(1)}(I)\) is generated by
      \[y(y-x)\del_y, \quad x(y-x)\del_x, \quad xy (\del_x + \del_y),
      \quad
      x\del_x + y \del_y.
      \]
      Since any product of distinct \(\delta_i\)'s is equal to
      \(P=xy(x-y)\), we have that for \(m \ge 2\)
      \begin{displaymath}
        \DO(I)^{(m)} = Sy(y-x)\del_y^m + Sx(y-x)\del_x^m + Sxy (\del_x +
        \del_y)^m + S\varepsilon_m + I\DO^{(m)}(S).
      \end{displaymath}
    \end{example}

    \end{subsubsection}

    \begin{subsubsection}{The case \(r \le n\)}
  For the generic arrangement \(\A_{n,r}\) with \(r \le n\), we can perform
  a linear change of variables so that \(p_i=x_i\). 
  Holm \cite[Paper III, Proposition 6.2]{Holm:Diff}
  showed the  following:
  \begin{proposition}\label{prop:rmnfree} 
  \(\DO^{(m)}(x_1\cdots x_r)\) is free with basis 
  \begin{math}
    \setsuchas{x_\alpha \del^\alpha}{\alpha \in \Nat^n, \, \tdeg{\alpha}=m},
  \end{math}
  where \(x_\alpha\) is the monic generator of 
  \(\sqrt{\langle x^{\tilde{\alpha}} \rangle}\), 
  the radical of \(\langle  x^{\tilde{\alpha}} \rangle\),
  and where 
  \begin{displaymath}
  \tilde{\alpha} =(\alpha_1,\dots,\alpha_r,0,\dots,0) \in \Nat^n.    
  \end{displaymath}
  \end{proposition}
  \end{subsubsection}

   \begin{subsubsection}{The case \(n=2\)}
    Holm \cite[Paper III, Prop 6.7]{Holm:Diff} notes that  an
    arrangement \(\A\) in \(\C^2\) is generic. Furthermore, he proves
    \begin{proposition}\label{prop:n2}
      Let \(\A\) be an arrangement in \(\C^2\), with defining
      polynomial \(p=p_1\cdots p_r \in \C[x_1,x_2]\),  and let \(m\)
      be a positive 
      integer. Put \(P_i = p/p_i\) for \(1 \le i \le r\), and define 
      \begin{displaymath}
        \delta_i = 
        \begin{cases}
          \del_2 & \text{ if } p_i = ax_1, a \in \C^*, \\
          \del_1 + a_i \del_2 & \text{ if } p_i = a(x_2-ax_1), a \in C^*
        \end{cases}
      \end{displaymath}
      Let 
      \begin{displaymath}
        \set{q_1,\dots,q_{\sbinom{n}{m}}} =
        \setsuchas{\del^\alpha}{\tdeg{\alpha}=m}. 
      \end{displaymath}
      Then \(\DO^{(m)}(\A)\) is free, minimally generated by 
      \begin{equation}
        \label{eq:n2gen}
        \begin{cases}
          \set{\varepsilon_m, P_1 \delta_1^m, \dots, P_r \delta_r^m} &
          \text{ if } 1 \le m \le r-2, \\
          \set{P_1\delta_1^{r-1}, \dots, P_r\delta_r^{r-1}} & \text{
            if } m=r-1, \\
          \set{P_1\delta_1^m, \dots, P_r\delta_r^m, pq_r,\dots,pq_m}& \text{
            if } m \ge r.
        \end{cases}
      \end{equation}
    \end{proposition}
  \end{subsubsection}

    \end{subsection}
  \end{section}

  \begin{section}{An exact sequence}
    If \(\A\) is an \(n\)-arrangement consisting of the
    hyperplanes \(H_1,\dots,H_r\), recall that the
    \emph{deleted} arrangement \(\A'\) is the arrangement in
    \(\C^n\) consisting of the hyperplanes \(H_2,\dots,H_{r}\). 
    The \emph{restricted} arrangement \(\A''\) is the arrangement
    \(H_2 \cap H_1, \dots, H_{r} \cap H_1 \subset H_1 \simeq
    \C^{n-1}\). Clearly, if \(\A\) is generic then so is \(\A'\) and
    \(\A''\), and we can write 
    \begin{displaymath}
      \A_{n,r}' = \A_{n,r-1}, \qquad \A_{n,r}''=\A_{n-1,r-1}.
    \end{displaymath}
      
    We can perform a change of coordinates so that the defining
    polynomial of \(H_1\) is \(x_n\). Then multiplication with
    \(x_n\) gives a \(S\)-module homomorphism of degree 1
    \begin{displaymath}
      \DO^{(m)}(\A_{n,r-1}) \xrightarrow{\cdot x_n} \rightarrow
      \DO^{(m)}(\A_{n,r}), 
    \end{displaymath}
    and the natural projection 
    \begin{displaymath}
      \pi: S=\C[x_1,\dots,x_n]   \to S' = \frac{S}{x_n} \simeq   
      \C[x_1,\dots,x_{n-1}]        
    \end{displaymath}
    induces a degree 0 map of graded vector spaces
    \begin{equation}\label{eq:pi}
      \DO^{(m)}(\A_{n,r}) 
      \rightarrow \DO^{(m)}(\A_{n-1,r-1})
    \end{equation}
    by replacing every occurrence of \(x_n\) by zero.
    The projection \(\pi\) can be used to give any \(S'\) module the
    structure of \(S\)-module (via extension of scalars), and so the
    map \eqref{eq:pi} becomes an \(S\)-module homomorphism.

    \begin{theorem}[\cite{Holm:Paper1}]
      If \(m,r,n\) are positive integers with \(n > 2\), and
      \(\A_{n,r}\) is a generic arrangement, then 
      there is a short exact sequence of graded \(S\)-modules
      \begin{equation}
        \label{eq:exact}
        0  
        \rightarrow \DO^{(m)}(\A_{n,r-1})(-1) 
        \rightarrow \DO^{(m)}(\A_{n,r}) 
        \rightarrow \DO^{(m)}(\A_{n-1,r-1})
        \rightarrow 0
      \end{equation}
    \end{theorem}

    \begin{cor}
      If \(m,r,n\) are positive integers with \(n > 2\) then 
      \begin{equation}
        \label{eq:hilbplus}
        \HI(\DO^{(m)}(\A_{n,r}))(b) = 
        b\HI(\DO^{(m)}(\A_{n,r-1}))(b) +
        \HI(\DO^{(m)}(\A_{n-1,r-1}))(b)
      \end{equation}
    \end{cor}

  \end{section}

  \begin{section}{The Hilbert and Poincar{\'e}-Betti series of
      \(\DO^{(m)}(\A)\) for generic \(\A\)}
    If \(\A\) and \(\A'\) are two generic \(n\)-arrangements with 
    \(\tdeg{\A} = \tdeg{\A'}=r\)
    then their Hilbert series and Poincar{\'e}-Betti series
    coincide. We let \(\A_{n,r}\) denote any generic \(n\)-arrangement
    consisting of \(r\) hyperplanes.

    \begin{subsection}{Derivations, the case \(m=1\)}
      Yuzvinsky \cite{Yuz:FreeRes} has given a minimal free resolution
      of \(\DO^{(1)}(\A_{n,r})\).
      \begin{theorem}[Yuzvinsky]
        Let \(r > n\ge 3\) and let \(\A_{n,r}\) be a generic
        \(n\)-arrangement with defining polynomial
        \(p=p_1\cdots p_r \in S=\C[x_1,\dots,x_n]\). Define
        \[\DO_0 = \setsuchas{\theta \in
          \DO^{(1)}(\A_{n,r})}{\theta(p_r)=0}.\]
        Then \( \DO^{(1)}(\A_{n,r}) = S\varepsilon_1 \oplus \DO_0\) as an
        \(S\)-module, and the minimal free resolution of \(\DO_0\) has
        length \(r-1\) and is \(r-n+1\)-linear.
        More precisely, the graded Betti numbers of \(\DO_0\) are given by 
      \begin{equation}
        \label{eq:der}
        \beta_{k,u} = 
        \begin{cases}
          \binom{r}{n-k} \binom{r-n+k-2}{k-1} & \text{ if } u+n-r-1=0
          \\
          0 & \text{ otherwise }
        \end{cases}
      \end{equation}
      \end{theorem}

      \begin{example}
      If \(n=3\), \(r=5\) then the Poincar{\'e}-Betti
      series of \(\DO_0\) is \(4b^3 + 2b^4t\), so the Poincar{\'e}-Betti
      series of \(\DO^{(1)}(\A_{3,5})\) is \(b+4b^3 + 2b^4t\).
      \end{example}

      We will give a compact formula for the Poincar{\'e}-Betti series
      of derivations, a formula that will give a hint as to what the
      Poincar{\'e}-Betti series of higher order differential operators
      might look like.

      \begin{lemma}
        Let \(r > n \ge 3\). Then 
        \begin{equation}
          \label{eq:pbder}
          \PB\left(\DO^{(1)}(\A_{n,r})\right)(b,t) = b + \left[ t^{n-r-1}(1+bt)^{r-1}
          \right] \hadamard (1-t)^{n-r}
        \end{equation}
      \end{lemma}

      \begin{example}
        To continue with the previous example, if \(n=3\), \(r=5\)
        then \eqref{eq:pbder} becomes 
        \begin{displaymath}
          \begin{split}
          \PB\left(\DO^{(1)}(\A_{3,5})\right) &= b + 
          \left[ t^{-3}(1+bt)^{4} \right] \hadamard (1-t)^{-2} \\
          & = b + \left[t^{-3} + 4bt^{-2} + 6 b^2t^{-1} + 4b^3 + b^4t
          \right] \hadamard (1 + 2t + 3t^2 + \dots) \\
          &= b + 4b^3 + 2b^4t
          \end{split}
        \end{displaymath}
      \end{example}
    \end{subsection}

    \begin{subsection}{The case \(n=2\)}
      Proposition~\ref{prop:n2} shows that when \(n=2\), the 
      Poincar{\'e}-Betti series of \(\DO^{(m)}(\A_{2,r})\) is given by 
       \begin{equation}
        \label{eq:n2pb}
        \PB(\DO^{(m)}(\A_{2,r}))(b,t) =
        \begin{cases}
          rb^{r-1} + (m-r+1)b^r & \text{ if } r \le 2 \text{ or } 
          m > r-2, \\
          b^m + mb^{r-1} & \text{ otherwise.}
        \end{cases}
      \end{equation}
      It follows that the Hilbert series is
      \begin{equation}
        \label{eq:n2hi}
        \HI(\DO^{(m)}(\A_{2,r}))(b,t) =
        \begin{cases}
          \frac{rb^{r-1} + (m-r+1)b^r}{(1-b)^2} & \text{ if } r \le 2
          \text{ or }  m > r-2, \\
          \frac{b^m + mb^{r-1}}{(1-b)^2} & \text{ otherwise.}
        \end{cases}
      \end{equation}
      Together with \eqref{eq:hilbplus}, the initial values \eqref{eq:n2hi}
      determines \(\HI(\DO^{(m)}(\A_{n,r}))(b)\) for
      all \(n \ge 2\).
 \end{subsection}

 \begin{subsection}{The case \(r \le n\)}
  Proposition~\ref{prop:rmnfree}  yields
  the following: 
  \begin{lemma}\label{lemma:pb}
    If \(r \le n\) then 
    \begin{equation}
      \label{eq:pbguess}
      \begin{split}
      \beta_{u,k}^{(m)}(\A_{n,r}) & = 
      \begin{cases}
        0 & \text{ if } u > 0 \\
        \binom{r}{k} \sbinom{n-r+k}{m-k} & \text{ if } u = 0
      \end{cases} \\
      \PB\left(\DO(\A_{n,r})\right)(a,b,t) &= \frac{(1-a+ab)^r}{(1-a)^n} \\
      \HI\left(\DO(\A_{n,r})\right)(a,b) &=
      \frac{(1-a+ab)^r}{(1-a)^n(1-b)^{-n}} \\ 
      \PB\left(\DO^{(m)}(\A_{n,r}\right)(b,t) &= [a^m]
      \frac{(1-a+ab)^r}{(1-a)^n} \\ 
      \HI\left(\DO^{(m)}(\A_{n,r}\right)(b) &= [a^m]
      \frac{(1-a+ab)^r}{(1-a)^n(1-b)^{-n}}  
      \end{split} 
    \end{equation}
  \end{lemma}
  \begin{proof}
    There are \(\binom{r}{k}\) different square-free \(\gamma \in
    \Nat^n\) of weight \(k\) such 
    that \[\gamma_{r+1}= \cdots = \gamma_{n}=0.\] 
    For a fixed such \(\gamma\), let \(g=g(\gamma)\) denote the
    number of \(\alpha \in \Nat^n\) with \(\bar{\alpha} = \gamma\). 
    Then 
    \[\alpha \mapsto \alpha -  \gamma\] gives a bijection between the
    set of multisets 
    \(\alpha\)
    on    \(\set{1,\dots,n}\) with weight \(k\) and with 
    \(\mathrm{supp}(\tilde{\alpha})= \gamma\), and
    the set of multisets on 
    \[\set{r+1,r+2,\dots,n} \cup \mathrm{supp}(\gamma)\] of weight \(m-k\). 
    So \(g\) is independent of \(\gamma\), 
    and is equal to \(\sbinom{n-r+k}{m-k}\).
    It follows that the
    number of minimal generators of \(\DO^{m}(\A)\) of polynomial degree
    \(k\) is given by \(\binom{r}{k} \sbinom{n-r+k}{m-k}\). 

    To prove the second identity, we take advantage of the fact that
    we may assume that the linear forms defining the arrangement are
    the monomials \(x_1,\dots,x_r\). For this particular arrangement,
    the modules of differentials will be multi-graded, by
    giving the operator
    \(x^\alpha \del^\beta\)
    the multi-degree \((\alpha,\beta)\). We will calculate
    Poincar{\'e}-betti series with respect to this fine grading, then
    specialize to get the desired series (which itself can not be
    given this fine grading).

    We start by calculating the generating
    function
    \begin{displaymath}
      G(a_1,\dots,a_n,b_1,\dots,b_n),
    \end{displaymath}
    where the sum is over  all minimal generators of \(\DO(\A_{n,r})\), and
    where where \(x^\alpha
    \del^\beta\) contributes \(a_1^{\beta_1} \cdots a_n^{\beta_n}
    b_1^{\alpha_1} \cdots b_n^{\alpha_n}\). Then as before
    we are looking for all pairs \((\alpha, \beta)\) with 
    \(\alpha \subset \set{0,\dots,r}\) and
    \(\mathrm{supp}(\tilde{\alpha})=\beta\). Hence 
    \begin{equation} \label{eq:pbgenrln}
      \begin{split}
        G(a_1,\dots,a_n,b_1,\dots,b_n) &= 
        \sum_{\beta \subseteq \set{1,2,\dots, r}} b^\beta 
        \left( \sum_{\mathrm{supp}(\gamma) \subset \beta} a^\gamma \right) 
        \left( 
          \sum_{ \mathrm{supp}(\theta) \subseteq \set{r+1,\dots,n}}
          a^\theta 
        \right) \\
        &= \sum_{\beta \subseteq \set{1,2,\dots, r}} a^\beta b^\beta 
        \prod_{i=r+1}^n (1-a_j)^{-1} \\
        &= \prod_{i=1}^r \frac{1-a_i+a_ib_i}{1-a_i} 
        \prod_{j=r+1}^n (1-a_j)^{-1}
      \end{split}
    \end{equation}
    Specializing gives 
    \begin{equation}
      \label{eq:spe}
      P(\A_{n,r})(a,b,t)=G(a,\dots,a,b,\dots,b)= \frac{(1-a+ab)^r}{(1-a)^n}
    \end{equation}
  \end{proof}
  
    \end{subsection}

    \begin{subsection}{The cases \(m=r-n+1\) and \(m > r-n+1\)}
      As we have already noted, it is known that \(\DO^{(m)}(\A_{n,r})\)
      is free when \(r \le n\) or when \(n \le 2\).
      Holm \cite[Paper III]{Holm:Diff} showed that 
      when \(r > n \ge 3\), \(m=r-n+1\), then \(\DO^{(m)}(\A_{n,r})\)
      is free module, minimally generated by \(\binom{r}{n-1}\)
      differential operators of order polynomial degree \(m\). 
      He conjectured that
      \(\DO^{(m)}(\A_{n,r})\) is free when 
      \(r > n \ge 3\) and \(m \ge r-n+1\).
      More precisely, it is reasonable to conjecture that \eqref{eq:pbguess}
      holds also for this range. This is certainly true for the Hilbert
      series:

      \begin{lemma}
        For \(r \le m+n-1\),
      \begin{equation}\label{eq:hfree}
       \HI(\DO^{(m)}(\A_{n,r}))(b) =  [a^m] \frac{(1-a+ab)^r}{(1-a)^n(1-b)^n} 
      \end{equation}
      \end{lemma}
      \begin{proof}
        This holds for \(n=2\) by \eqref{eq:n2hi}. Furthermore, if
        \(r \le m+n-1\) then \(r-1 \le m+n-1\) and \(r-1 \le
        m+n-1-1\), so
      the assertion follows by induction, since
      \begin{displaymath}
        [a^m] \frac{(1-a+ab)^r}{(1-a)^n(1-b)^n} =
        [a^m] \frac{(1-a+ab)^{r-1}}{(1-a)^n(1-b)^n} +
        [a^m] \frac{(1-a+ab)^{r-1}}{(1-a)^{n-1}(1-b)^{n-1}}.        
      \end{displaymath}
      \end{proof}
    \end{subsection}

    \begin{subsection}{The case \(r \ge m+n\)}
      Holm \cite[Paper III]{Holm:Diff} showed that 
      when \(r > n \ge 3\), \(r \ge m+n\), then
      \(\DO^{(m)}(\A_{n,r})\)
      is \emph{not} a free module.
      This is therefore the ``interesting range''. We'll eventually
      formulate a conjecture regarding the Poincar{\'e}-Betti series
      of \(\DO^{(m)}(\A_{n,r})\) for \(m,n,r\) in this range.
      
      We'll simplify the problem slightly by identifying a direct
      summand of these modules.
      Recall the notations of Theorem~\ref{thm:gengen}, so that 
      \(D=\set{P_1\delta_1, \dots, P_t\delta_t}\) is a certain subset
      of \(\DO^{(1)}(\A_{n,r})\) with the property that 
      \begin{displaymath}
        \setsuchas{P_\alpha \delta^\alpha}{\alpha \in \Nat^t, \,
          \tdeg{\alpha}=m} \cup \set{\varepsilon_m}
      \end{displaymath}
      generates \(\DO^{(m)}(\A_{n,r})\), where \(P_\alpha\) is the product
      of those \(p_i\) such that some \(\delta_j\) with \(\alpha_j
      \neq 0\) does not annihilate \(p_i\).
      
      Let \(\Xi^{(m)}(\A_{n,r})\) be the module generated by 
      \begin{displaymath}
        \setsuchas{P_\alpha \delta^\alpha}{\alpha \in \Nat^t, \,
          \tdeg{\alpha}=m}.
      \end{displaymath}
      Then Holm's result can be stated as
      \begin{equation}
        \label{eq:summydum}
      \DO^{(m)}(\A_{n,r}) = \Xi^{(m)}(\A) + S\varepsilon_m        
      \end{equation}

      Holm showed \cite[Paper I, Lemma 5.27]{Holm:Diff} that 
      for \(r \le n\), \(\varepsilon_m \in \Xi^{(m)}(\A)\). 
      Furthermore, we have:
      \begin{lemma}\label{lemma:summand}
        Suppose that \(r > n \ge 3\), \(r \ge m+n\). Then
        \(\varepsilon_m \not \in \Xi^{(m)}(\A_{n,r})\), so 
        \begin{equation}
          \label{eq:exsum}
          \DO^{(m)}(\A_{n,r}) = \Xi^{(m)}(\A) \oplus S\varepsilon_m
        \end{equation}
        Hence, 
        \begin{equation}
          \label{eq:pbplus}
          \PB\left(\DO^{(m)}(\A_{n,r})\right)(b,t) = b^m +
          \PB\left(\Xi^{(m)}(\A_{n,r})\right)(b,t) 
        \end{equation}
      \end{lemma}
      \begin{proof}
        All  \(P_i\) have degree \(r-n+1\), so all \(P_\alpha\) have
        degrees \(\ge r-n+1 > m\), hence all \(P_\alpha
        \delta^\alpha\) have polynomial degrees \(> m\). But
        \(\varepsilon_m\) have polynomial degree \(m\), hence can not
        be expressed as an \(S\)-linear combination of the \(P_\alpha
        \delta^\alpha\)'s. This shows that \(\varepsilon_m \not \in
        \Xi^{(m)}(\A)\), which together with \eqref{eq:summydum}
        yields \eqref{eq:exsum}.
      \end{proof}

      Below, we tabulate
      \(\PB(\Xi^{(m)}(\A_{n,r}))(b,t)\) for \(m=2,3,4\) (these
      Poincar{\'e}-Betti series were
      calculated using Macaulay 2 \cite{MACAULAY2} and the method described in the
      beginning).

      For \(m=2\), the Poincar{\'e}-Betti series are as follows:

        \begin{center}
          \begin{Small}
      \begin{tabular}{|c|ccc|}
        \hline
        r & n=3 & n=4 &n=5 \\ \hline &&&\\
        5 &   $2 b^{4} {t}+7 b^{3}$   & -    & -  \\
        6 &  $4 b^{5} {t}+9 b^{4}$ & $2 b^{5} t^{2}+9 b^{4}
        {t}+16 b^{3}$ & -   \\ 
        7 &  $6 b^{6} {t}+11 b^{5}$ & $6 b^{6} t^{2}+22 b^{5}
        {t}+25 b^{4}$ &  $2 b^{6} t^{3}+11 b^{5} t^{2}+25 b^{4}
        {t}+30 b^{3}$ \\  
        8 & $8 b^{7} {t}+13 b^{6}$ & $12 b^{7} t^{2}+39 b^{6}
        {t}+36 b^{5}$ & $8 b^{7} t^{3}+39 b^{6} t^{2}+72 b^{5}
        {t}+55 b^{4}$ \\ 
        9 & $10 b^{8} {t}+15 b^{7}$  & $20 b^{8} t^{2}+60 b^{7}
        {t}+49 b^{6}$ & \\ 
        10 &$12 b^{9} {t}+17 b^{8}$ &&\\  \hline
      \end{tabular}
          \end{Small}
        \end{center}
      
      We conjecture that 
      \begin{multline}
        \label{eq:pdm2}
        \PB(\Xi^{(2)}(\A_{n,r})(b,t) = 
        \left[
          t^{n-r-1}(1+2bt)(1+bt)^{r-2}
          \right] \hadamard (1-t)^{n-r+1}
      \end{multline}

      For \(m=3\), we get
        \begin{center}
          \begin{Small}
            \begin{tabular}{|c|ccc|}
              \hline
              r & n=3 & n=4 & n=5 \\ \hline &&&\\
              6 & $3 b^{5} {t}+12 b^{4}$ & - & - \\
              7 & $6 b^{6} {t}+15 b^{5}$ & $3 b^{6} t^{2}+15
              b^{5} {t}+31 b^{4}$ & - \\
              8 & $9 b^{7} {t}+18 b^{6}$ &
              $9 b^{7} t^{2}+36 b^{6} {t}+46 b^{5}$ &
               $3 b^{7} t^{3}+18 b^{6} t^{2}+46 b^{5} {t}+65
               b^{4}$ \\
               9 & $12 b^{8} {t}+21 b^{7}$ 
               & $18 b^{8} t^{2}+63 b^{7} {t}+64 b^{6}$
               & \( 111b^5 + 128tb^6 + 63t^2b^7 + 12t^3b^8\) \\ \hline
            \end{tabular}
          \end{Small}
        \end{center}
        We conjecture that 
      \begin{multline}
        \label{eq:pdm3}
        \PB(\Xi^{(3)}(\A_{n,r}))(b,t) = \\
        \left[
          t^{n-r-1}(1+3bt + 3b^2t^2)(1+bt)^{r-3}
          \right] \hadamard (1-t)^{n-r+2}
      \end{multline}

      For \(m=4\) we get 
      \begin{center}
        \begin{small}
          \begin{tabular}{|c|cc|}
            \hline
            r & n=3 & n= 4 \\ \hline &&\\
            7 & \(18b^5 +4tb^6\) & - \\
            8 & \(52b^5 +22tb^6 + 4t^2b^7\)&
            $8 b^{7} {t}+22 b^{6}$ \\ \hline
          \end{tabular}
        \end{small}
      \end{center}
        We conjecture that 
      \begin{multline}
        \label{eq:pdm4}
        \PB(\Xi^{(4)}(\A_{n,r}))(b,t) = \\
        \left[
          t^{n-r-1}(1+4bt + 6b^2t^2 +4b^3t^3)(1+bt)^{r-4}
          \right] \hadamard (1-t)^{n-r+3}
      \end{multline}

      Based on these computations, we make the following conjecture,
      which by Yuzvinsky's result is true for derivations, i.e. when
      \(m=1\). 
      \begin{conjecture}
        Suppose that \(3 \le n\) and \(r \ge m+n\). Let \(\A_{n,r}\)
        be a generic \(n\)-arrangement with \(\tdeg{\A_{n,r}}=r\). Then
      \begin{multline}
        \label{eq:pdmm}
        \PB\left(\Xi^{(m)}(\A_{n,r})\right) =  \\
        \left\{ t^{-r +n - 1} \bigl((1+bt)^m -(bt)^m \bigr) 
        \bigr(1+bt\bigr)^{r-m} \right\} 
      \hadamard (1-t)^{m - r+n - 1}
      \end{multline}
      and hence
      \begin{multline}
        \PB\left(\DO^{(m)}(\A_{n,r})\right) =  \\
        b^m+
      \left\{ t^{-r +n - 1} \bigl((1+bt)^m -(bt)^m \bigr) 
        \bigr(1+bt\bigr)^{r-m} \right\} 
       \hadamard (1-t)^{m - r+n - 1}
      \end{multline}
      \end{conjecture}
      
      Note that this conjecture implies that
      \begin{enumerate}[(i)]
      \item The homological dimension of the \(S\)-module
       \(\Xi^{(m)}(\A_{n,r})\) (and of \(\DO^{(m)}(\A_{n,r})\)) is \(n-2\).
      \item The minimal free resolution  of
        \(\Xi^{(m)}(\A_{n,r}))\) is \(r-n\)-linear.
      \item  \(\Xi^{(m)}(\A_{n,r})\) is minimally generated by
        \begin{displaymath}
          \binom{r}{r-n+1} - \binom{r-m}{r-n+1-m} =
          \binom{r}{n-1} -    \binom{r-m}{n-1} 
        \end{displaymath}
        differential operators of polynomial
        degree \(r-n+1\); \(\DO^{(m)}(\A_{n,r})\) is minimally
        generated by these differential operators and \(\varepsilon_m\).
      \end{enumerate}



    We now indicate a possible way of proving {eq:pdmm}.
    The short exact sequence \eqref{eq:exact}, together with
    Lemma~\ref{lemma:summand}, gives a short exact sequence       
    \begin{equation}
      \label{eq:exactred}
      0  
      \rightarrow \Xi^{(m)}(\A_{n,r-1})(-1) 
      \rightarrow \Xi^{(m)}(\A_{n,r}) 
      \rightarrow \Xi^{(m)}(\A_{n-1,r-1})
      \rightarrow 0
    \end{equation}
    which  gives rise to a long
    exact sequence in homology (we have omitted the shifts which are
    necessary to make the morphisms below homogeneous of degree zero)
    \begin{multline}
      \label{eq:longexact} 
      \cdots \rightarrow
      \Tor^1(\DO^{(m)}(\A_{n,r-1}), \C) \rightarrow 
       \Tor^1(\DO^{(m)}(\A_{n,r}), \C) \rightarrow \\
       \rightarrow \Tor^1(\DO^{(m)}(\A_{n-1,r-1}), \C) \xrightarrow{\delta_1} 
       \Xi^{(m)}(\A_{n,r-1}) \tensor \C  \rightarrow\\
      \rightarrow \Xi^{(m)}(\A_{n,r}) \tensor \C
      \rightarrow \Xi^{(m)}(\A_{n-1,r-1}) \tensor \C
      \rightarrow 0
    \end{multline}
    which controls the ``deviation'' 
    \begin{displaymath}
      q(m,n,r)= \PB(\Xi^{(m)}(\A_{n,r})) -
      b\PB(\Xi^{(m)}(\A_{n,r-1})) -
      (1+bt)\PB(\Xi^{(m)}(\A_{n-1,r-1})).
    \end{displaymath}
    If all the connecting homomorphisms \(\delta_i\) are zero, then
    so is this deviation, and the Poincar{\'e}-Betti series can be
    computed recursively using deletion-restriction. 
    In the ``interesting range'' \(3 \le n < r\), \(r \ge m+n\),
    assuming the conjectured formula \eqref{eq:pdmm}, 
    and using Lemma~\ref{lemma:had} below, we do get that \(q(m,n,r) =
    0\). This indicates (but does not prove) that the  connecting
    homomorphisms are zero. 
    Conversely, a proof of the vanishing of all connecting
    homomorphisms  would also prove \eqref{eq:pdmm}.

    \begin{lemma}\label{lemma:had}
      For \(3 \le n < r\), \(r \ge m+n\), it holds that
      \begin{multline}
        \label{eq:cons1}
        \left\{ t^{-r +n - 1} \bigl((1+bt)^m -(bt)^m \bigr) 
          \bigr(1+bt\bigr)^{r-m} \right\} \hadamard (1-t)^{m - r+n - 1} \\
        - 
        b\left\{ t^{-r +n } \bigl((1+bt)^m -(bt)^m \bigr) 
          \bigr(1+bt\bigr)^{r-1-m} \right\} \hadamard (1-t)^{m - r+n
        } 
        \\
        -
        (1+bt) \left\{ t^{-r +n-1 } \bigl((1+bt)^m -(bt)^m \bigr) 
          \bigr(1+bt\bigr)^{r-1-m} \right\} \hadamard (1-t)^{m -
          r+n-1 }  = 0
      \end{multline}
    \end{lemma}
    \begin{proof}[Proof of Lemma~\ref{lemma:had}]
      Put \(k=r-n\) and 
      \begin{displaymath}
        \begin{split}
          U(r,k,m) & = \left[ \bigl( (1+bt)^m -(bt)^m \bigr) 
            \bigr(1+bt\bigl)^{r-m} \right] \\
          &= t^{-k-1}(1+bt)^r - t^{m-k-1}(1+bt)^{r-m}
        \end{split}
      \end{displaymath}
      Then for \(\ell > 0\) we have that 
        \begin{align*}
          [x^\ell] && U(r,k,m) &= & b^{\ell+k+1} & \left[
            \binom{r}{\ell+k+1} - \binom{r-m}{\ell-m+k+1} \right] \\
          [x^\ell]  && b \, U(r-1,k-1,m) &=& b^{\ell+k+1} & \left[
            \binom{r-1}{\ell+k} - \binom{r-1-m}{\ell-m+k} \right] \\
          [x^\ell] && U(r-1,k,m) &=& b^{\ell+k+1} &\left[
            \binom{r-1}{\ell+k+1} - \binom{r-1-m}{\ell-m+k+1}
          \right] \\
          [x^\ell] && bt \, U(r-1,k,m) &=& b^{\ell+k+1} & \left[
            \binom{r-1}{\ell+k} - \binom{r-1-m}{\ell-m+k} \right] 
        \end{align*}
      Since 
      \begin{displaymath}
        (1-t)^{m-r+n-1} = \sum_{\ell=0}^\infty \sbinom{m-r+n-1}{\ell}
      \end{displaymath}
      the equation \eqref{eq:cons1} is equivalent to the identity
      \begin{multline}
        \label{eq:binid}
        \left[ 
          \binom{r}{\ell+k+1} - \binom{r-m}{\ell-m+k+1}
        \right] \sbinom{m-k+1}{\ell} \\
        - \left[
          \binom{r-1}{\ell+k} - \binom{r-1-m}{\ell-m+k} 
        \right] \sbinom{m-k}{\ell} \\
        - \left[
          \binom{r-1}{\ell+k+1} - \binom{r-1-m}{\ell-m+k+1}
        \right] \sbinom{m-k+1}{\ell} \\
        - \left[
          \binom{r-1}{\ell+k} - \binom{r-1-m}{\ell-m+k} 
        \right] \sbinom{m-k+1}{\ell-1} =0,
      \end{multline}
      which can be verified algorithmically\footnote{We used the
        \texttt{simplify(expr,GAMMA)} command of the computer
        algebra system Maple \cite{MAPLE}.}.
      
      The case \(\ell=0\) is dealt with similarly.
    \end{proof}

  \end{subsection}

\end{section}

\bibliographystyle{plain}
\bibliography{journals,articles}

\begin{thebibliography}{10}

\bibitem{LouAD:GB}
William~W. Adams and Philippe Loustaunau.
\newblock {\em An introduction to {G}r\"obner bases}, volume~3 of {\em Graduate
  Studies in Mathematics}.
\newblock American Mathematical Society, Providence, RI, 1994.

\bibitem{Bjork:DiffOp}
J.-E. Bj{\"o}rk.
\newblock {\em Rings of differential operators}, volume~21 of {\em
  North-Holland Mathematical Library}.
\newblock North-Holland Publishing Co., Amsterdam, 1979.

\bibitem{MAPLE}
Bruce~W. Char, Keith~O. Geddes, Gaston~H. Gonnet, Benton~L. Leong, Michael~B.
  Monagan, and Stephen~M. Watt.
\newblock {\em {Maple V Library reference manual.}}
\newblock Springer-Verlag, 1991.

\bibitem{Coutinho:Primer}
S.~C. Coutinho.
\newblock {\em A primer of algebraic {$D$}-modules}, volume~33 of {\em London
  Mathematical Society Student Texts}.
\newblock Cambridge University Press, Cambridge, 1995.

\bibitem{MACAULAY2}
Daniel~R. Grayson and Michael~E. Stillman.
\newblock Macaulay 2.
\newblock Computer algebra program, available at
  \url{http://www.math.uiuc.edu/Macaulay2/}.

\bibitem{Holm:Diff}
P{\"a}r Holm.
\newblock {\em Differential {O}perators on {A}rrangements of {H}yperplanes}.
\newblock PhD thesis, {S}tockholm {U}niversity, 2002.

\bibitem{Holm:Paper1}
P{\"a}r Holm.
\newblock Differential operators on hyperplane arrangements.
\newblock {\em Communications in {A}lgebra}, 32(6):2177 -- 2201, 2004.
\newblock DOI: 10.1081/AGB-120037213.

\bibitem{BasicDer}
Tadeusz J{\'o}zefiak and Bruce~E. Sagan.
\newblock Basic derivations for subarrangements of {C}oxeter arrangements.
\newblock {\em J. Algebraic Combin.}, 2(3):291--320, 1993.

\bibitem{GeomFreeArr}
Joseph P.~S. Kung.
\newblock A geometric condition for a hyperplane arrangement to be free.
\newblock {\em Adv. Math.}, 135(2):303--329, 1998.

\bibitem{OrlikTerao:ArrHyp}
Peter Orlik and Hiroaki Terao.
\newblock {\em Arrangements of hyperplanes}, volume 300 of {\em Grundlehren der
  Mathematischen Wissenschaften [Fundamental Principles of Mathematical
  Sciences]}.
\newblock Springer-Verlag, Berlin, 1992.

\bibitem{RoseTerao}
Lauren~L. Rose and Hiroaki Terao.
\newblock A free resolution of the module of logarithmic forms of a generic
  arrangement.
\newblock {\em Journal of {A}lgebra}, 136:376--400, 1991.

\bibitem{Dmod}
Harry Tsai and Anton Leykin.
\newblock D-modules.m2.
\newblock Macaulay 2 package, available at
  \url{http://www2.math.uic.edu/~leykin/Research/Dmodules/Dmods.html}, 2002.

\bibitem{Yuz:FreeRes}
Sergey Yuzvinsky.
\newblock A free resolution of the module of derivations for generic
  arrangements.
\newblock {\em Journal of {A}lgebra}, 136:432--436, 1991.

\end{thebibliography}

\end{document}